\documentclass[11pt,final]{amsart}

\usepackage[
  a4paper,
  margin=2.5cm
]{geometry}
\usepackage[T1]{fontenc}
\usepackage[utf8]{inputenc}
\usepackage{lmodern}
\usepackage{microtype}
\usepackage{mathtools}
\usepackage{amssymb}
\usepackage{mathrsfs}
\usepackage{bm}
\usepackage{enumitem}
\usepackage{booktabs}
\usepackage{xcolor}
\usepackage{hyperref}
\usepackage{aliascnt}
\usepackage[nameinlink,capitalize,noabbrev]{cleveref}

\hypersetup{
colorlinks=true,
linkcolor=black
}

\numberwithin{equation}{section}

\newcommand{\R}{\mathbb{R}}

\newcommand{\dd}{\,\mathrm{d}}

\newtheorem{theorem}{Theorem}[section]
\newaliascnt{proposition}{theorem}
\newtheorem{proposition}[proposition]{Proposition}
\aliascntresetthe{proposition}
\newaliascnt{lemma}{theorem}
\newtheorem{lemma}[lemma]{Lemma}
\aliascntresetthe{lemma}
\newaliascnt{corollary}{theorem}

\aliascntresetthe{corollary}
\theoremstyle{definition}
\newaliascnt{definition}{theorem}

\aliascntresetthe{definition}
\theoremstyle{remark}
\newaliascnt{remark}{theorem}
\newtheorem{remark}[remark]{Remark}
\aliascntresetthe{remark}
\title{
Antisymmetry Breaking for Critical Poincar\'e--Sobolev Extremals
in the Three Dimensional Ball
}

\author{Zhang Yitian}

\date{}

\begin{document}

\begin{abstract}
We study extremal functions for Poincar\'e--Sobolev inequalities in the unit ball in $\R^3$.
At the critical exponent $q=6$, we prove that no minimizer is antisymmetric with respect to any 
hyperplane through the origin. We obtain a positive gap between the general and antisymmetric optimal constants. 
Continuity of these constants yields antisymmetry breaking for exponents sufficiently close to $6$. 
We introduce two critical indices to describe the onset of antisymmetry breaking and prove that global minimizers 
are antisymmetric at both indices.
\end{abstract}

\maketitle

\section{Introduction}
\label{sec:introduction}

Let $B=B_1(0)\subset\R^3$. Consider the optimal Poincar\'e--Sobolev constants
\[
\Lambda_q
\coloneqq
\inf_{\substack{0\ne u\in H^1(B)\\ \int_Bu\,\dd x=0}}
\dfrac{\int_B|\nabla u|^2\,\dd x}
{\left(\int_B|u|^q\,\dd x\right)^{2/q}},
\qquad 2\leqslant q\leqslant 6.
\]
For $q=2$, the minimizers are first nonconstant Neumann eigenfunctions and are hyperplane antisymmetric. 
The question is whether this property persists when $q>2$.

Gir\~ao and Weth~\cite{GiraoWeth2006} proved that every minimizer is axially symmetric and strictly decreasing 
in the polar angle from a suitable axis. They also proved antisymmetry with respect to the hyperplane through 
the origin perpendicular to this axis when $q$ is close to $2$. In dimension two, they established antisymmetry breaking 
for $q$ large enough. In dimensions $N \geqslant 3$, they left open whether minimizers fail to be antisymmetric; 
see~\cite[Remark~1.4(v)]{GiraoWeth2006}. Related symmetry questions were studied for noncoercive functionals by 
Brock, Croce, Guib\'e, and Mercaldo~\cite{BrockCroceGuibeMercaldo2020} and for weighted Sobolev quotients by Brock, Chiacchio, Croce,
and Mercaldo~\cite{BrockEtAl2025}. Both works establish antisymmetry breaking in dimension two.

The one dimensional version is well understood. Building on earlier work, Buslaev, Kondratev and Nazarov~\cite{Nazarov2002,BuslaevKondratevNazarov1998}
completed the classification of mean zero extremals. For the quotient on $[0,1]$, minimizers are antisymmetric with
respect to $\frac{1}{2}$ for $2\leqslant q\leqslant 6$ and are not  antisymmetric for $q>6$. Ghisi, Gobbino, and Rovellini~\cite{GhisiGobbinoRovellini2018} 
later gave the sharp symmetry criterion for generalized Wirtinger inequalities.

We answer this question in dimension three near the critical exponent $q=6$. Let $\Lambda_q^{\mathrm{odd}}$ denote the 
optimal constant restricted to functions antisymmetric with respect to a fixed hyperplane through the origin. Let $S_3$ be 
the sharp Euclidean Sobolev constant~\cite{Aubin1976,Talenti1976}. Theorem~\ref{thm:critical_antisymmetry_breaking} establishes
\begin{equation}
\Lambda_6
<2^{-2/3}S_3
<\frac{2}{3}S_3
\leqslant\Lambda_6^{\mathrm{odd}}.
\end{equation}
The upper bound for $\Lambda_6$ and the existence of critical minimizers are due to~\cite[Proposition~1.2]{GiraoWeth2006}. 
We get the lower bound for $\Lambda_6^{\mathrm{odd}}$, which gives a positive gap between the two optimal constants. 
Consequently, no critical minimizer is antisymmetric with respect to any hyperplane through the origin. 

The proof combines Kelvin transform, a mixed Steklov inequality, and the sharp Euclidean Sobolev inequality. We restrict an 
antisymmetric function to a half ball and extend it to the half space by Kelvin transform. The ground state representation 
bounds the spherical boundary term in the Kelvin energy identity by the Dirichlet energy. The required ground state and 
Kelvin identities follow from \cite{PinchoverVersano2023} and \cite{BrezisLieb1985,CostabelDallaRivaDaugeMusolino2025}, respectively.

We also obtain two critical indices in $(2,6)$ and prove antisymmetry breaking for exponents sufficiently close to $6$. 
The argument combines continuity of the optimal constants, antisymmetry near $2$, and the positive gap. At both indices, 
the optimal constants coincide, and compactness yields antisymmetric global minimizers.

Section~\ref{sec:notation} introduces the notations.
Section~\ref{sec:symmetry_breaking} proves antisymmetry breaking when $q=6$. 
Section~\ref{sec:critical_indices} establishes continuity of the optimal constants and studies the critical indices.

\section{Notations}
\label{sec:notation}

We write
\[
\|u\|_q
\coloneqq
\left(\int_B|u|^q\,\dd x\right)^{1/q},
\qquad
\mathcal{H}
\coloneqq
\left\{
u\in H^1(B):
\int_Bu\,\dd x=0
\right\}.
\]
By the Poincaré inequality,
\[
\|u\|_{\mathcal{H}}
\coloneqq
\|\nabla u\|_{L^2(B)}
\]
is a well defined norm on $\mathcal{H}$.

For $2\leqslant q\leqslant 6$, define
\[
\mathcal{Q}_q(u)
\coloneqq
\dfrac{\int_B|\nabla u|^2\,\dd x}
{
 \left(\int_B|u|^q\,\dd x\right)^{2/q}},
\qquad
0\ne u\in\mathcal{H},
\]
and
\[
\Lambda_q
=
\Lambda_q(B)
\coloneqq
\inf_{0\ne u\in\mathcal{H}}\mathcal{Q}_q(u)>0.
\]
The normalized minimizing set is
\[
\mathcal{M}_q
\coloneqq
\left\{
u\in\mathcal{H}:
\|u\|_q=1,\ 
\int_B|\nabla u|^2\,\dd x=\Lambda_q
\right\}.
\]
It is nonempty for $q<6$ by compactness and for $q=6$ by \cite[Proposition~1.2]{GiraoWeth2006}.

For $e\in\mathbb S^2$, let
\[
\sigma_e(x)
\coloneqq
x-2(x\cdot e)e,
\]
and define
\[
H^1_{\mathrm{odd}}(B;e)
\coloneqq
\left\{
u\in H^1(B):
u\circ\sigma_e=-u
\right\}
\subset\mathcal{H}.
\]
The restricted optimal constant is
\[
\Lambda_q^{\mathrm{odd}}
=
\Lambda_q^{\mathrm{odd}}(B;e)
\coloneqq
\inf_{0\ne u\in H^1_{\mathrm{odd}}(B;e)}
\mathcal{Q}_q(u)
\geqslant\Lambda_q.
\]
Since $B$ and $\mathcal{Q}_q$ are rotationally invariant, $\Lambda_q^{\mathrm{odd}}(B;e)$ is independent of $e$.
Thus, whenever a particular antisymmetry plane is considered, we may assume without loss of generality that
\[
e=e_3\coloneqq(0,0,1).
\]
In this case,
\[
\sigma_{e_3}(x_1,x_2,x_3)=(x_1,x_2,-x_3).
\]

The normalized odd minimizing set is
\[
\mathcal{M}_q^{\mathrm{odd}}(e)
\coloneqq
\left\{
u\in H^1_{\mathrm{odd}}(B;e):
\|u\|_q=1,\ 
\int_B|\nabla u|^2\,\dd x
=
\Lambda_q^{\mathrm{odd}}
\right\},
\]
which is nonempty for $q<6$.

We further set
\[
\begin{aligned}
\Delta(q)
&\coloneqq
\Lambda_q^{\mathrm{odd}}-\Lambda_q,
\\
\mathcal{A}
&\coloneqq
\{q\in(2,6]:\Delta(q)>0\},
\\
q_c
&\coloneqq
\inf\mathcal{A}.
\\
\mathcal{C}
&\coloneqq
\left\{
p\in(2,6]:
\Delta(q)>0
,\, \forall q\in[p,6]
\right\},
\\
q^c
&\coloneqq
\inf\mathcal{C}.
\end{aligned}
\]

For functions axially symmetric about the $x_3$-axis, we use
\[
r=|x|,
\qquad
\theta
=
\arccos\left(\frac{x_3}{|x|}\right),
\]
and write $u=u(r,\theta)$. Reflection across $\left\{x\in\mathbb{R}^3\colon x_3=0\right\}$ corresponds to
\[
(r,\theta)\longmapsto(r,\pi-\theta),
\]
so that antisymmetry is equivalent to
\[
u(r,\pi-\theta)=-u(r,\theta).
\]

For the Kelvin transform argument, we use
\[
\R^3_+
=
\{x\in\R^3:x_3>0\},
\qquad
B^+
=
B\cap\R^3_+,
\]
and
\[
D
=
B\cap \left\{x\in\mathbb{R}^3\colon x_3=0\right\},
\qquad
\Sigma^+
=
\partial B\cap\{x_3>0\}.
\]
We use the mixed boundary space
\[
H^1_D(B^+)
\coloneqq
\left\{
v\in H^1(B^+):
\operatorname{Tr}_Dv=0
\right\}
\]
and, for $v\in H^1_D(B^+)$, set
\[
E(v)
\coloneqq
\int_{B^+}|\nabla v|^2\,\dd x,
\qquad
L(v)
\coloneqq
\int_{B^+}|v|^6\,\dd x,
\qquad
T(v)
\coloneqq
\int_{\Sigma^+}v^2\,\dd\sigma.
\]
The Kelvin extension of $v$ is defined by
\[
(\mathcal{K}v)(x)
\coloneqq
\begin{dcases}
v(x),
&x\in B^+,
\\
|x|^{-1}v\left(\frac{x}{|x|^2}\right),
&x\in\R^3_+\setminus\overline B.
\end{dcases}
\]

Finally, let
\[
\dot H^1(\R^3)
\coloneqq
\overline{C_c^\infty(\R^3)}^{\,\|\nabla\cdot\|_2},
\qquad
\dot H^1_0(\R^3_+)
\coloneqq
\overline{C_c^\infty(\R^3_+)}^{\,\|\nabla\cdot\|_2},
\]
and the sharp Euclidean Sobolev constant is
\[
S_3
\coloneqq
\inf_{0\ne w\in\dot H^1(\R^3)}
\dfrac{\int_{\R^3}|\nabla w|^2\,\dd x}
{\left(\int_{\R^3}|w|^6\,\dd x\right)^{1/3}}.
\]

\section{Antisymmetry Breaking}
\label{sec:symmetry_breaking}

The following identity is a direct specialization of the mixed boundary ground state transform in \cite[the proof of Theorem~6.1]{PinchoverVersano2023}.

\begin{proposition}
\label{prop:mixed_Steklov_representation}
For every $v\in H_D^1(B^+)$, the weighted gradient
\[
x_3\nabla\left(\frac{v}{x_3}\right)
\]
is well defined as an element of $L^2(B^+;\R^3)$, and
\begin{equation}
\label{eq:ground_state_general}
\int_{B^+}|\nabla v|^2\,\dd x
-
\int_{\Sigma^+}v^2\,\dd\sigma
=
\int_{B^+}
x_3^2
\left|
\nabla\left(\frac{v}{x_3}\right)
\right|^2
\dd x.
\end{equation}
Consequently,
\[
\int_{\Sigma^+}v^2\,\dd\sigma
\leqslant
\int_{B^+}|\nabla v|^2\,\dd x,
\]
with equality if and only if
\[
v(x)=c\,x_3
\]
for some $c\in\R$.
\end{proposition}

\begin{proof}
In our setting, the bilinear form in \cite[Definition~2.7]{PinchoverVersano2023} is
\[
\mathcal{B}(f,g)
\coloneqq
\int_{B^+}\nabla f\cdot\nabla g\,\dd x
-
\int_{\Sigma^+}fg\,\dd\sigma.
\]
First let $v$ be smooth on $\overline{B^+}$ and vanish near $\overline D$, and set
\[
u(x)=x_3,
\qquad
\phi=\frac{v}{x_3}.
\]
Denote $\nu$ the outward unit normal to $\partial B$. Since
\[
-\Delta u=0\quad\text{in }B^+,
\qquad
u=0\quad\text{on }D,
\qquad
\partial_\nu u=u\quad\text{on }\Sigma^+,
\]
the ground state representation in \cite[the proof of Theorem~6.1]{PinchoverVersano2023} gives
\[
\mathcal{B}(v,v)
=
\int_{B^+}u^2|\nabla\phi|^2\,\dd x
=
\int_{B^+}x_3^2
\left|\nabla\left(\frac{v}{x_3}\right)\right|^2\,\dd x.
\]
This is \eqref{eq:ground_state_general}. By density and continuity of the trace map, it extends to every $v\in H_D^1(B^+)$, 
and simultaneously defines $x_3\nabla(v/x_3)$ as an element of $L^2(B^+;\R^3)$.

Finally, the equality 
\[
\int_{\Sigma^+}v^2\,\dd\sigma 
= 
\int_{B^+}|\nabla v|^2\,\dd x
\]
holds if and only if
\[
\nabla\left(\frac{v}{x_3}\right)=0
\quad\text{a.e. in }B^+,
\]
which, since $B^+$ is connected, is equivalent to $v=cx_3$.
\end{proof}

\begin{theorem}
\label{thm:critical_antisymmetry_breaking}
No minimizer of $\Lambda_6$ is antisymmetric with respect to any hyperplane through the origin.

More precisely, for every $e\in\mathbb S^2$ and every nonzero $u\in H^1_{\mathrm{odd}}(B;e)$,
\begin{equation}
\label{eq:quantitative_odd_gap}
\mathcal{Q}_6(u)-\Lambda_6
>
\left(
\frac{2}{3}-2^{-2/3}
\right)S_3.
\end{equation}
In particular, no minimizing sequence contained in any fixed hyperplane antisymmetric class can be a minimizing sequence for $\Lambda_6$.
\end{theorem}

\begin{proof}
By rotational invariance, it suffices to consider
\[
0\ne u\in H^1_{\mathrm{odd}}(B;e_3).
\]
Set $v=u|_{B^+}$. Then $v\in H_D^1(B^+)$, and antisymmetry gives
\begin{equation}
\label{eq:odd_half_ball_identities}
\int_B|\nabla u|^2\,\dd x
=
2E(v),
\qquad
\int_B|u|^6\,\dd x
=
2L(v).
\end{equation}

Let $V=\mathcal{K}v$. Formula~(3.21) in \cite{BrezisLieb1985}, applied to the zero extension of $v$ to $B$, yields
\begin{equation}
\label{eq:kelvin_L6}
\int_{\R^3_+}|V|^6\,\dd x=2L(v).
\end{equation}
Moreover, \cite[Lemma~A.1 and formula~(A.3)]{CostabelDallaRivaDaugeMusolino2025}, with
$n=3$, $\Gamma=\R^3_+$, and $\rho=1$, gives the exterior Kelvin energy identity. Since the interior and exterior pieces have matching traces on $\Sigma^+$ and zero trace on $\{x_3=0\}$, they glue to a function $V\in\dot H^1_0(\R^3_+)$ satisfying
\begin{equation}
\label{eq:kelvin_energy}
\int_{\R^3_+}|\nabla V|^2\,\dd x
=
2E(v)+T(v).
\end{equation}

Applying the sharp Sobolev inequality to the zero extension of $V$ to $\R^3$, we obtain
\[
S_3\left(2L(v)\right)^{1/3}
=
S_3\left(\int_{\R^3_+}|V|^6\,\dd x\right)^{1/3}
\leqslant
\int_{\R^3_+}|\nabla V|^2\,\dd x
=
2E(v)+T(v).
\]
By Proposition~\ref{prop:mixed_Steklov_representation}, $T(v)\leqslant E(v)$. Hence
\[
S_3\left(2L(v)\right)^{1/3}
\leqslant
3E(v),
\]
and therefore
\[
\frac{E(v)}{L(v)^{1/3}}
\geqslant
\frac{2^{1/3}}{3}S_3.
\]
Using \eqref{eq:odd_half_ball_identities}, we conclude that
\begin{equation}
\label{eq:odd_critical_lower_bound}
\mathcal{Q}_6(u)
=
2^{2/3}\frac{E(v)}{L(v)^{1/3}}
\geqslant
\frac{2}{3}S_3.
\end{equation}

On the other hand, \cite[Proposition~1.2]{GiraoWeth2006} yields
\[
\Lambda_6<2^{-2/3}S_3.
\]
Together with \eqref{eq:odd_critical_lower_bound}, this gives
\[
\mathcal{Q}_6(u)-\Lambda_6
>
\left(
\frac{2}{3}-2^{-2/3}
\right)S_3,
\]
proving \eqref{eq:quantitative_odd_gap}.
\end{proof}

\section{Critical Indices for Antisymmetry Breaking}
\label{sec:critical_indices}

We first establish the continuity of the two optimal constants.

\begin{lemma}
\label{lem:continuity_optimal_constants}
The maps
\[
q\longmapsto\Lambda_q,
\qquad
q\longmapsto\Lambda_q^{\mathrm{odd}}
\]
are continuous on $[2,6]$. Consequently,
\[
q\longmapsto\Delta(q)
\]
is continuous on $[2,6]$.
\end{lemma}

\begin{proof}
We first prove the continuity of $\Lambda_q$ at $q_0\in[2,6)$. Let $q_j\to q_0$, and fix $r\in(q_0,6)$. For all sufficiently large $j$, we have $q_j<r$, so we may choose $u_j\in\mathcal{M}_{q_j}$.

For every fixed $0\ne\varphi\in\mathcal{H}$,
\[
\|\nabla u_j\|_2^2
=
\Lambda_{q_j}
\leqslant
\mathcal{Q}_{q_j}(\varphi)
=
\frac{\|\nabla\varphi\|_2^2}{\|\varphi\|_{q_j}^2}.
\]
Since $\varphi\in L^r(B)$ and $q_j\to q_0<r$, the dominated convergence theorem gives
\[
\|\varphi\|_{q_j}\longrightarrow\|\varphi\|_{q_0}.
\]
Thus $(\nabla u_j)$ is bounded in $L^2(B)$, and the Poincaré inequality shows that $(u_j)$ is bounded in $H^1(B)$.

Passing to a further subsequence if necessary, we may assume that
\[
u_j\rightharpoonup u
\quad\text{in }H^1(B),
\qquad
u_j\to u
\quad\text{in }L^r(B).
\]
Since $q_j<r$,
\[
\|u_j-u\|_{q_j}
\leqslant
|B|^{1/q_j-1/r}\|u_j-u\|_r
\longrightarrow0.
\]
Moreover, $u\in L^r(B)$ and $q_j\to q_0$, so
\[
\|u\|_{q_j}\longrightarrow\|u\|_{q_0}.
\]
Consequently,
\[
\left|
\|u_j\|_{q_j}-\|u\|_{q_0}
\right|
\leqslant
\|u_j-u\|_{q_j}
+
\left|
\|u\|_{q_j}-\|u\|_{q_0}
\right|
\longrightarrow0.
\]
Since $\|u_j\|_{q_j}=1$, it follows that $\|u\|_{q_0}=1$. Moreover, $\mathcal{H}$ is weakly closed in $H^1(B)$, 
and hence $u\in\mathcal{H}\setminus\{0\}$. Weak lower semicontinuity now yields
\[
\Lambda_{q_0}
\leqslant
\int_B|\nabla u|^2\,\dd x
\leqslant
\liminf_{j\to\infty}\Lambda_{q_j}.
\]

Conversely, taking $w\in\mathcal{M}_{q_0}$, we obtain
\[
\limsup_{j\to\infty}\Lambda_{q_j}
\leqslant
\lim_{j\to\infty}\mathcal{Q}_{q_j}(w)
=
\mathcal{Q}_{q_0}(w)
=
\Lambda_{q_0}.
\]
Therefore,
\[
\Lambda_{q_j}\longrightarrow\Lambda_{q_0}.
\]

It remains to prove continuity at $q=6$. For $2\leqslant q<6$, H\"older's inequality gives
\[
\|u\|_q
\leqslant
|B|^{1/q-1/6}\|u\|_6.
\]
It follows that
\[
\Lambda_q
\geqslant
|B|^{1/3-2/q}\Lambda_6.
\]
Hence
\[
\liminf_{q\uparrow6}\Lambda_q\geqslant\Lambda_6.
\]

Conversely, for every $\varepsilon>0$, choose $0\ne w\in\mathcal{H}$ such that
\[
\mathcal{Q}_6(w)\leqslant\Lambda_6+\varepsilon.
\]
Then
\[
\limsup_{q\uparrow6}\Lambda_q
\leqslant
\lim_{q\uparrow6}\mathcal{Q}_q(w)
=
\mathcal{Q}_6(w)
\leqslant
\Lambda_6+\varepsilon.
\]
Letting $\varepsilon\to0$, we conclude that
\[
\Lambda_q\longrightarrow\Lambda_6
\qquad\text{as }q\uparrow6.
\]

The same argument, with $\mathcal{H}$ replaced by the closed subspace $H^1_{\mathrm{odd}}(B;e)$, proves the continuity of $\Lambda_q^{\mathrm{odd}}$. 
The continuity of $\Delta$ follows immediately.
\end{proof}

We can now locate the two critical indices.

\begin{theorem}
\label{thm:existence_critical_indices}
The sets $\mathcal{A}$ and $\mathcal{C}$ are nonempty, and
\[
2<q_c\leqslant q^c<6.
\]
Moreover,
\begin{align}
\Lambda_q
&=
\Lambda_q^{\mathrm{odd}},
\quad
\forall 2\leqslant q\leqslant q_c,
\label{eq:equality_below_qc}
\\
\Lambda_q
&<
\Lambda_q^{\mathrm{odd}},
\quad
\forall q^c<q\leqslant 6,
\label{eq:breaking_above_qc}
\end{align}
and
\[
\Delta(q_c)=\Delta(q^c)=0.
\]
In particular,  $\forall e\in\mathbb S^2$,
\[
\varnothing\ne
\mathcal{M}_{q_c}^{\mathrm{odd}}(e)
\subset
\mathcal{M}_{q_c},
\qquad
\varnothing\ne
\mathcal{M}_{q^c}^{\mathrm{odd}}(e)
\subset
\mathcal{M}_{q^c}.
\]
\end{theorem}

\begin{proof}
By \cite[Theorem~1.3(d)]{GiraoWeth2006} and rotational invariance, there exists $\varepsilon_0>0$ such that
\[
\Delta(q)=0,
\quad
\forall 2\leqslant q\leqslant 2+\varepsilon_0.
\]
On the other hand, Theorem~\ref{thm:critical_antisymmetry_breaking} gives
\[
\Delta(6)
\geqslant
\left(
\frac{2}{3}-2^{-2/3}
\right)S_3
>0.
\]
By Lemma~\ref{lem:continuity_optimal_constants}, there exists $\varepsilon_1>0$ such that
\[
\Delta(q)>0
\qquad
\text{for }6-\varepsilon_1\leqslant q\leqslant 6.
\]
Consequently, $\mathcal{A}$ and $\mathcal{C}$ are nonempty and
\[
2<q_c\leqslant q^c<6.
\]

By the definition of $q_c$,
\[
\Delta(q)=0
\qquad
\text{for }2\leqslant q<q_c.
\]
Continuity then gives $\Delta(q_c)=0$, proving \eqref{eq:equality_below_qc}.

The set $\mathcal{C}$ is upward closed. Hence every $q\in(q^c,6]$ belongs to $\mathcal{C}$, and therefore $\Delta(q)>0$, proving \eqref{eq:breaking_above_qc}. 
If $\Delta(q^c)>0$, continuity would imply that
\[
\Delta(q)>0,
\qquad
\forall q\in[q^c-\delta,6]
\]
for some $\delta>0$. This would imply that $q^c-\delta\in\mathcal{C}$, which is a contradiction. Thus $\Delta(q^c)=0$.

Finally, since $q_c,q^c<6$, compactness gives
\[
\mathcal{M}_{q_c}^{\mathrm{odd}}(e)\ne\varnothing,
\qquad
\mathcal{M}_{q^c}^{\mathrm{odd}}(e)\ne\varnothing.
\]
Since
\[
\Lambda_{q_c}^{\mathrm{odd}}=\Lambda_{q_c},
\qquad
\Lambda_{q^c}^{\mathrm{odd}}=\Lambda_{q^c},
\]
we obtain
\[
\mathcal{M}_{q_c}^{\mathrm{odd}}(e)\subset\mathcal{M}_{q_c},
\qquad
\mathcal{M}_{q^c}^{\mathrm{odd}}(e)\subset\mathcal{M}_{q^c}.
\]
\end{proof}

\begin{remark}
The theorem does not assert that $q_c=q^c$. This equality is equivalent to the absence of antisymmetry restoration after antisymmetry breaking.
\end{remark}

\section*{Declaration of AI use}
During the preparation of this manuscript, the author used OpenAI Chatgpt Plus without Codex to assist with drafting and revising mathematical arguments, 
checking calculations and the consistency of proofs, improving the organization and language of the manuscript, and editing the \LaTeX{} source. 
The authors take full responsibility for the mathematical results, proofs, references, and final content of the manuscript.

\end{document}